\magnification1200
\input pst-plot

\font\tener=eurm10
\font\sevener=eurm7
\font\fiveer=eurm5
\textfont1=\tener
\scriptfont1=\sevener
\scriptscriptfont1=\fiveer

\mathcode`0="7130
\mathcode`1="7131
\mathcode`2="7132
\mathcode`3="7133
\mathcode`4="7134
\mathcode`5="7135
\mathcode`6="7136
\mathcode`7="7137
\mathcode`8="7138
\mathcode`9="7139

\font\tenfrak=eufm10
\font\sevenfrak=eufm7
\font\fivefrak=eufm5
\newfam\frakfam
\textfont\frakfam=\tenfrak
\scriptfont\frakfam=\sevenfrak
\scriptscriptfont\frakfam=\fivefrak

\font\tenam=msam10
\font\sevenam=msam7
\font\fiveam=msam5
\newfam\amssymfam
\textfont\amssymfam=\tenam
\scriptfont\amssymfam=\sevenam
\scriptscriptfont\amssymfam=\fiveam

\count255=\the\amssymfam
\multiply\count255 by"100
\advance\count255 by"2000
\advance\count255 by"03\mathchardef\square   =\count255
\def\qed{{\hfill$\square$}}

\font\tentxexa=txexa
\font\seventxexa=txexa scaled 700
\font\fivetxexa=txexa scaled 500
\newfam\txexafam
\textfont\txexafam=\tentxexa
\scriptfont\txexafam=\seventxexa
\scriptscriptfont\txexafam=\fivetxexa

\count255=\the\txexafam
\multiply\count255 by"100
\advance\count255 by"1000
\advance\count255 by"10\mathchardef\bigtimes=\count255

\font\sevenrm=tir at 7.6pt

\rm

\font\sevenbf=tib at 7.6pt

\nopagenumbers
\def\title{On Normal Subgroups of Coxeter Groups}
\font\sevenss=cmss7
\headline{\sevenss \ifnum\pageno>1
\ifodd\pageno\title\ \hrulefill\ \the\pageno
\else\the\pageno\ \hrulefill\ S.~R.~Gal
\fi\fi}
\centerline{\bf \title}
\centerline{\bf Generated by Standard Parabolic Subgroups}
\footnote{}{2000 {\it Mathematics Subject Classification: }20F55}
\smallskip
\centerline{\'Swiatos\l aw R. Gal%
\footnote{$^\star$}{Partially supported by a KBN grant 2 P03A 017 25.}}
\centerline{Wroc\l aw University}
\font\tt=cmtt10
\centerline{\tt http://www.math.uni.wroc.pl/\~{}sgal/papers/nscg.dvi}
\bigskip

{
\narrower\narrower\smallskip\everymath{\scriptstyle}\sevenrm

\noindent{\sevenbf Abstract :}
We discuss one construction of nonstandard subgroups in
the category of Coxeter groups.

\noindent Two formulae for the growth series of such a subgroups are given.

\noindent As an application we construct a flag simple convex polytope,
whose f-polynomial has non-real roots.

\smallskip}   

\parindent0pt
\parskip=\smallskipamount
\newcount\secno\secno=0
\newcount\thmno\thmno=0

\def\section #1\par{\vskip 0pt plus.1\vsize\penalty-250
\vskip 0pt plus-.1\vsize\bigskip\vskip\parskip
\global\advance\secno by 1\global\thmno=0
{\tenfrak\the\secno.} {\bf #1}\par
\smallskip
}

\def\tag{\the\secno.\the\thmno}
\def\mktag{{\global\advance\thmno by 1}\tag}

\beginsection Introduction

The central object of this paper is the growth series of a Coxeter group.
Many geometric features of such a group (or any group in general)
reflect in properties of the growth series.

We describe in detail the normal closure of a standard parabolic in a
Coxeter group $W$. It is again a Coxeter group and its Coxeter presentation
is given explicitly.

We give two formulae for the growth series of such a normal subgroup.
The first one is given as a specialization of a multi-variable version
of the growth series of $W$.

The second formula works when the normal subgroup is right-angled (this can
be easily checked by analysis of the Dynkin diagram of $W$).
The formula is based on  counting
special subdiagrams of the Dynkin diagram of $W$.

As an application we construct a flag simple convex polytope,
whose f-polynomial has non-real roots.

The author would like to thank Tadeusz Januszkiewicz
for useful discussions,
Pawe{\l} Goldstein for his extensive help with the
final version of this manuscript and
Mike Davis for a preliminary version of [D2].

\section Preliminaries on Coxeter Groups

\proclaim Definition. A {\rm Coxeter system} $(W,S)$ is a group $W$
together with a set of generators $S$, and presentation
$$W=\langle S\vert (st)^{m_{st}}=1\hbox{\rm\ for all }s,t\in S\rangle,$$
where $m_{ss}=1$ (i.e. all the generators have order two)
and $m_{st}=m_{ts}\in\{2,3,\ldots,\infty\}$ if $s\neq t$.
One reads $(st)^\infty=1$ as no relation between $s$ and $t$ imposed.
The matrix $m$ is called the {\rm Coxeter matrix} of\/ $W$.
We usually ignore $S$ and call $W$ a {\rm Coxeter group}. 
The subgroup of\/ $W$ generated by $T\subset S$ is denoted $W_T$. 
Such a subgroup is called a {\rm standard parabolic}
or {\rm parabolic} for short.

Traditionally a Coxeter matrix is depicted in a following decorated graph
called {\it Dynkin diagram}. Its vertex set is $S$ and any two distinct
vertices are joined by an edge with label $m_{st}$. By convention
{\parindent=.3in\parskip=0pt
\item{$\bullet$}we omit an edge if $m_{st}=2$,
\item{$\bullet$}we omit a label if $m_{st}=3$,
\item{$\bullet$}we draw double edge instead that labeled by $4$ and
\item{$\bullet$}we draw dashed edge instead that labeled by $\infty$.
}

We will use the following convention for the Cayley graph: the group acts on
the left, therefore $\{g_1,g_2\}$ is an edge provided $g_1^{-1}g_2\in S$.
We say that an edge $\{g_1,g_2\}$ is {\it marked} by the generator
$g_1^{-1}g_2$.

\proclaim Theorem \mktag\ [D2]. Let $W$ be a group
generated by a set of involutions $S$. Let $C_W$ be the Cayley
graph of $(W,S)$. Then $(W,S)$ is a Coxeter system if and only if
for any element $s\in S$ the fixpoint set
$(C_W)^s$ separates $C_W$ (i.e. $C_W-(C_W)^s$ has two components
being interchanged by $s$).

\edef\davisthm{Theorem \tag}

{\it Proof: } It is a part of Theorem 2.3.3 in [D2]
(cf. [D2, Definitions 2.2.1 and 2.2.10]).\qed

\proclaim Definition. A conjugate of a generator from $S$ is called
{\rm a reflection}. The set of all reflections is denoted by $R(W)$.

\proclaim Corollary \mktag. Reflections are intrinsically
defined by a Coxeter system. Precisely, if $T\subset S$ then
$R(W_S)\cap W_T=R(W_T)$.

\edef\reflcor{Corollary \tag}

{\it Proof: } Assume that $r\in R(W_S)\cap W_T$.
Since the action of $r$ is conjugate to that of some generator,
the fixpoint set of $r$ also separates $C_W$. Then $r$ and the neutral
element $e$ are in different components of $C_{W_S}-(C_{W_S})^r$, and thus
in different components of $C_{W_T}-(C_{W_T})^r$. Take any path joining
$e$ and $r$ in $C_{W_T}$. Such a path has to intersect the separating
fixpoint set $(C_{W_T})^r$ in a midpoint of some edge $\{w,wt\}$.
Thus $rw=wt$ and, what follows, $r=wtw^{-1}\in R(W_T)$.\qed

\proclaim Definition. The {\rm length} of an element $w$ of any group
{\rm with respect to a given generating set} $S$ is denoted by $\ell(w)$
and is equal to the minimal length of a word in $S$ representing $w$.
The word is called {\rm reduced} if its length is equal to the length
of the element it represents.

\proclaim Lemma \mktag\ [T]. If $(W,S)$ is a Coxeter system
than any word can be reduced  (into any reduced word representing
the same element of\/ $W$)
in a sequence of moves of the following form:
{\parindent=.3in\parskip=0pt
\item{(a)} removing a subword of a form $ss$, or
\item{(b)} replacing an alternating subword $st\ldots$ of length $m_{st}$
by an alternating subword $ts\ldots$ of the same length.
}

\edef\titslemma{Lemma \tag}

{\it Note: } The move (a) is not symmetric, i.e. we do not need creation of
a pair $ss$. In particular any two reduced words representing the same
element may be joined by moves of the form (b) only. 

\section Normal closure of a standard parabolic

If $m_{st}=2n+1$ is odd then $s$ and $t$ are conjugates since
the relation reads then $$(st)^ns(st)^{-n}=t.$$

If one wants to find a normal closure of a subgroup $W_T$
it is convenient to incorporate in $T$ all the generators conjugate
to those already in $T$. This clarifies the assumption of the following


\proclaim Proposition \mktag. Assume that $T\subset S$ is such that
for any $t\in T$ and $s\in S-T$ the exponent $m_{st}$ is even
(or $\infty$). Then
{\parskip=0pt\parindent=.2in
\item{(1)} there exists a well defined homomorphism $\phi_T\colon W\to W_{S-T}$
that is an identity on $W_{S-T}$ and sends $W_T$ to the identity element,
\item{(2)} the normal closure $\overline{W_T}$ of\/ $W_T$
is the  kernel of $\phi_T$,
\item{(3)} $\overline{W_T}$ is generated by $S_T$, a set of all the conjugates
of $T$ with respect to $W_{S-T}$,
\item{(4)} the Cayley graph $C_{\overline{W_T}}$ of
$(\overline{W_T},S_T)$ is made from the Cayley graph $C_W$ of\/ $W$
by collapsing all the edges marked by the generators not in $T$,
\item{(5)} $(\overline{W_T},S_T)$ is a Coxeter system.
}

\edef\mainthm{Proposition \tag}

{\it Proof :} Let $t\in T$ and $s\not\in T$, then the relation
$st\ldots=ts\ldots$, after substituting the identity for $s$,
becomes trivial exactly when $m_{st}$ is even. Thus (1) follows.
This also proves that the relation ``$s$ and $t$ are conjugate'' is
the equivalence closure of the relation ``$m_{st}$ is odd''.

Any element $w$ of $W$ may be written as follows
$$\eqalign{&w_0t_1w_1\ldots t_nw_n=\cr
&\hphantom{=}
\left(w_0t_1w_0^{-1}\right)\left(w_0w_1t_2(w_0w_1)^{-1}\right)\ldots
\left(w_0\ldots w_{n-1}t_n(w_0\ldots w_{n-1})^{-1}\right)
\left(w_0\ldots w_n\right)\cr}$$
where $t_i\in T$ and $w_i\in W_{S-T}$. If $e=\phi_T(w)=w_0\ldots w_n$
then $w$ is expressible in $S_T$, therefore (2) and (3).

Let $g_1u_1$, $g_2u_2$, where $g_i\in\overline{W_T}$ and $u_i\in W_{S-T}$,
be the endpoints of an edge corresponding to some generator in $T$,
i. e. $(g_1u_1)^{-1}g_2u_2\in T$.
In particular $e=\phi_T(u_1^{-1}g_1^{-1}g_2u_2)=\phi_T(u_1^{-1}u_2)$.
Therefore $u_1=u_2$. Furthermore
$$g_1^{-1}g_2=u_1\left((g_1u_1)^{-1}(g_2u_1)\right)u_1^{-1}\in S_T.$$
On the contrary, if $g_1^{-1}g_2=utu^{-1}$, then the edge $\{g_1,g_2\}$
in $C_{\overline{W_T}}$ is covered by the edge $\{g_1u,g_2u\}$
in $C_W$. This proves (4).

Let $t\in S_T$.
By \davisthm\ we need to show that the fixpoint set of $t$ separates
$C_{\overline{W_T}}$. Let $W^+$ and $W^-$ denote the two components of
$C_W-(C_W)^t$. Two vertices, one from $W^+$ and the other from $W^-$
are joined by an edge precisely when they differ by $t$ on the right.
Since $t$ is not a conjugate of any element of $W_{S-T}$ all the cosets
of $W_{S-T}$ have to lie in $W^+$ or $W^-$. Therefore the fixpoint set
of $t$ also separates $C_{\overline{W_T}}$. Thus (5).\qed

{\it Remark \mktag:} If $W_{S-T}$ is not finite,
then $\overline{W_T}$ has in general an infinite number of generators.

\edef\finrem{Remark \tag}

\section Coxeter matrix of $\overline{W_T}$

\def\supp{\mathop{\rm supp}}

\proclaim Definition. Define the support $\supp w$ of an element
$w$ of a Coxeter group to be the set of generators that appear
in its reduced presentation. By \titslemma\ this does not depend
on the reduced presentation.

\proclaim Definition. Let $T_1,T_2\subset S$. We write $T_1\perp T_2$
if $m_{t_1t_2}=2$ for any $t_1\in T_1$ and $t_2\in T_1$. If $T\subset S$
we define $T^\perp=\{s\in S\colon m_{st}=2\hbox{\rm{\ for all }} t\in T\}$.  

Fix $t_1,t_2\in T$.
Let $w^*$ be the shortest
element in the double coset $W_{t_1^\perp-T}wW_{t_2^\perp-T}$.
Such an element is unique [B, Ch. IV, \S1, Ex. 3].

\proclaim Corollary \mktag. Let $\tau_1=w_1t_1w_1^{-1}$
and $\tau_2=w_2t_2w_2^{-1}$ be two generators of $\overline{W_T}$.
Then $\tau_1=\tau_2$ if and only if
$$t_1=t_2\quad\hbox{\sl and}\quad(w_1^{-1}w_2)^*=1.$$
If $\tau_1\neq\tau_2$ then we have
$$m_{\tau_1\tau_2}=
\cases{m_{t_1t_2}&if $t_1\neq t_2$ and $(w_1^{-1}w_2)^*=1$,\cr
m_{st_1}\over 2&if $t_1=t_2$ and $(w_1^{-1}w_2)^*=s\in S-T$,\cr
\infty& otherwise.\cr}$$

\edef\coxmat{Corollary \tag}

{\it Proof: }
Set $w=(w_1^{-1}w_2)^*$, so $$w_1^{-1}w_2=\eta_1w\eta_2$$
for some $\eta_i\in W_{t_i^\perp-T}$. In particular $t_i$
commutes with $\eta_i$ 

Let $N$ be a natural number such that
$$(\tau_1\tau_2)^N=1.\leqno (\mktag)$$
\edef\relation{($\tag$)}%
Then \relation\ reads
$$\eqalign{1=\eta_1^{-1}w_1^{-1}w_1\eta_1
&=\eta_1^{-1}w_1^{-1}
\left(w_1t_1(w_1^{-1}w_2)t_2\left(w_1^{-1}w_2\right)^{-1}\ldots\right)w_1\eta_1\cr
&=\eta_1^{-1}
\left(t_1(\eta_1w\eta_2)t_2(\eta_1w\eta_2)^{-1}\ldots\right)
\eta_1\cr
&=(\eta_1^{-1}t_1\eta_1)w(\eta_2t_2\eta_2^{-1})w^{-1}\ldots\cr
&=t_1wt_2w^{-1}\ldots\cr}$$

\titslemma\ implies that $w$ has length at most one, otherwise
no cancellation is possible since any reduced form of $w$ does not start
with a letter that commutes with $t_1$ nor
ends with a letter that commutes with $t_2$.

The case $\ell(w)=0$ is obvious. If $w$ is a generator, then it
commutes neither with $t_1$ nor with $t_2$. The only possibility
to use \titslemma\
is described by the claim.\qed

\section Multi-variable growth series.

The length function on a Coxeter group may be refined
in the following way. Take a set of generators $S_0\subset S$
with the property that whenever $s$ belongs to $S_0$, then
all the generators conjugate to $s$ belong to $S_0$.

\titslemma\ implies the following

\proclaim Corollary. The number (counted with multiplicities) of generators
from $S_0$ appearing in a reduced presentation of $w$, and denoted
$\ell_{S_0}(w)$, does not depend on the presentation.

Note that $\ell_{S_0}$ satisfies the triangle inequality, i.e.
$\ell_{S_0}(w_1w_2)\leq\ell_{S_0}(w_1)+\ell_{S_0}(w_2)$.

Let $\varphi\colon S\to I$ be any partition of the generating set
into the sets with the above property.
We call such a partition {\it allowable}.
For each $w\in W$
define a monomial in indeterminate $I$:
$w_\varphi({\bf x})=\prod_{x\in I}x^{\ell_{\varphi^{-1}(x)}(w)}$, where
${\bf x}=(x)_{x\in I}$.

For any subset $A$ of $W$, a formal power series
$A_\varphi({\bf x})=\sum_{w\in A}w_\varphi({\bf x})$ is called
a multi-variable growth series of $A\subset W$.

The above definition needs some finiteness assumption on $(W,S)$.
It is enough to assume that $S$ is finite (if not, the
coefficients of the growth series may be infinite). Since we
also want to study the growth series of $\overline{W_T}$,
we want to know that $S_T$ is finite. Therefore, due to \finrem,
we will assume that $W_{S-T}$ is finite.

The notation may look a little ambiguous.
However, note that $W_\varphi$ does not depend on whether
$W$ is a standard  parabolic subgroup of some bigger Coxeter group or not.
Similarly $\left(\overline{W_T}\right)_\varphi$ may denote two things:
either a growth series of $\overline{W_T}$ or that of $\overline{W_T}
\subset W$. In fact these two series coincide, as shown by the following

\proclaim Proposition \mktag. Assume that $T\subset S$ is as in \mainthm\ 
and $\varphi^{-1}(x_0)=S-T$.
There is a well defined prolongation $\varphi_T$ of $\varphi$ to $S_T$
by requiring that $\varphi(wtw^{-1})=\varphi(t)$ for any $t\in T$
and $w\in W_{S-T}$. It is obviously allowable. Then
$$(\overline{W_T})_{\varphi_T}({\bf x})=
(\overline{W_T})_\varphi({\bf x})_{|x_0=1},$$
where the left series is a growth series of $\overline{W_T}$
and on the right that of $\overline{W_T}\subset W$.

\edef\wwprop{Proposition \tag}

{\it Note: }The right hand side is well defined,
i.e. there is only a finite number of nonzero coefficients at $Jx_0^m$
for given monomial $J$ in variables different from $x_0.$ 

If $w_\varphi({\bf x})=Jx_0^m$ then $w=w_0t_{j_1}w_1\ldots w_{n-1}t_{j_n}w_n$,
where $\prod\varphi(t_{j_k})=J$ and $w_k\in W_{S-T}$.
There is only a finite number of possibilities for $t_{j_k}$ and $w_k$
(notice that we have assumed that $W_{S_T}$ is finite).

{\it Proof  of \wwprop:}
Take a reduced word in $S$ representing an element from $\overline{W_T}$
and write it as $$w_0t_1w_1\ldots t_nw_n\leqno(\mktag)$$
where $t_i\in T$ and $w_j$ are subwords in $S-T$.

\edef\wred{(\tag)}

If the following word in $S_T$ (representing the same element)
$$\left(w_0t_1w_0^{-1}\right)\left(w_0w_1t_2(w_0w_1)^{-1}\right)\ldots
\left(w_0\ldots w_{n-1}t_n(w_0\ldots w_{n-1})^{-1}\right)\leqno(\mktag)$$
\edef\wtred{(\tag)}%
were not reduced, one would do some cancellations in \wtred, expand the result
to obtain the word in $S$ with
smaller $\ell_{\varphi^{-1}(x)}$ for some $x\neq x_0$ than that of \wred.
This would show that \wred\ was not reduced.

Finally \wred\ and \wtred\ define the same monomial. Thus the claim. \qed 

\proclaim Theorem \mktag.
Under the assumptions of the preceding proposition
$$W_\varphi({\bf x})_{|x_0=1}=(\# W_{S-T})\cdot
\left(\overline{W_T}\right)_{\varphi_T}({\bf x}).$$

\edef\specthm{Theorem \tag}

{\it Proof :}
Any element $w$ of $W$ decomposes uniquely as $w=w_Tw_0$ where,
$w_T\in\overline{W_T}$ and $w_0\in W_{S-T}$. We need to show that
$\ell_{\varphi^{-1}(x)}(w)=\ell_{\varphi^{-1}(x)}(w_T)$ for any $x\neq x_0$.

Since $\ell_{\varphi^{-1}(x)}(w_0)=0$, by the triangle inequality we have
$$\ell_{\varphi^{-1}(x)}(w_T)\geq
\ell_{\varphi^{-1}(x)}(w_Tw_0)\geq
\ell_{\varphi^{-1}(x)}(w_Tw_0w_0^{-1})=\ell_{\varphi^{-1}(x)}(w_T).$$
\qed

{\it Note: } Two generators of $\overline{W_T}$ conjugate
in $W$ may not be conjugate in $\overline{W_T}$, thus not every multi variable
growth series of $\overline{W_T}$ is specialization of that
of $W$.

\proclaim Corollary. Assume that $T\subset S$ is as in \mainthm.
Let $\varphi\colon S\to\{x_0,x\}$ be such that $\varphi(t)=x$ if
and only if $t\in T$. Then the ordinary growth series of\/ $W_T$ and
$\overline{W_T}$ are computed as follows:
$$\eqalign{W_T(x)&=W_\varphi(0,x),\cr
\overline{W_T}(x)&={W_\varphi(1,x)\over W_\varphi(1,0)},\cr}$$
where $W_\varphi(x_0,x)$ is a multi-variable growth series of\/ $W$
associated to the allowable partition $\varphi$.

{\it Proof: } $W_\varphi(1,0)=W_{S-T}(1)=\# W_{S-T}$.

\section Right-angled Coxeter groups and flag complexes

\proclaim Definition. {\rm A nerve} $N_W$ of a Coxeter system $(W,S)$
is a simplicial complex
consisting of $T\subset S$ such that the subgroup $W_T$ is finite.

{\it Remark : } Some authors define
a nerve as a baricentric subdivision of the above.

\proclaim Definition \mktag.
A Coxeter group is said to be {\rm right-angled} if $m_{st}\in\{1,2,\infty\}$
for all $s,t\in S$.

If $W$ is a right-angled Coxeter group, then $W_T$ is
finite if and only if $W_{T'}$ is finite for any two-element subset $T'$
of $T$, therefore $N_W$ is a {\it flag completion} of its one-skeleton:
$T$ is a face of $N_W$ if and only if its one skeleton is contained in $N_W$.
Such a simplicial complex is called {\it a flag complex}.

On the other hand, let $\Gamma$ be any graph. Let $S$ denote the set
of vertices of $\Gamma$. Declare
$$m_{st}=\cases{1&if $s=t$,\cr
2&if there is an edge joining $s$ and $t$,\cr
\infty&otherwise.\cr}$$
Then $N_W$ is the flag completion of $\Gamma$.

Let $X$ be any simplicial complex.
For a given function $\varphi\colon S\to I$ on the set of vertices
one defines an f-polynomial of $X$ by the formula
$$f_{X,\varphi}({\bf x})\colon=
\sum_{\sigma\in X}\prod_{s\in\sigma}\varphi(s).
\leqno(\mktag)$$

\proclaim Proposition \mktag\ [S, Prop.~26]. Let $W$ be an arbitrary
Coxeter group and $\varphi$ an allowable partition. Then
$W_\varphi({\bf x})$ is a series of
a rational function. Moreover, if $W$ is infinite,
then $${1\over W_\varphi({\bf x}^{-1})}=
\sum_{T\subset S}{(-1)^{\#T}\over\left(W_T\right)_{\varphi_{|T}}({\bf x})},$$
where $T$ runs over subsets of $S$ such that $W_T$ is finite, and
${\bf x}^{-1}=(x^{-1})_{x\in I}$.

\edef\serre{Proposition \tag}

\proclaim Corollary \mktag. Assume that $W$ is a right angled Coxeter group.
Since no two generators are conjugated, any function
$\varphi\colon S\to I$ is an allowable partition and
$$f_{N_W,\varphi}\left({-1\over1+{\bf x}}\right)
={1\over W_{\varphi}({\bf x}^{-1})}. \leqno(\mktag)$$

\edef\fw{(\tag)}

{\it Proof: }
If $W_T$ is finite, then
$\left(W_T\right)_{\varphi}({\bf x})=
\prod_{t\in T}(1+\varphi(t))$, thus \fw\ follows form \serre.\qed

\section Computation of f-polynomial

\coxmat\ implies the following

\proclaim Corollary \mktag. $\overline{W_T}$ is a right-angled Coxeter group
if and only if
{\parindent=.3in\parskip=0pt
\item{\rm(1)} if $s\not\in T$ and $t\in T$ then $m_{st}\in\{2,4,\infty\}$,
and
\item{\rm(2)} if $t,t'\in T$ then $m_{tt'}\in\{1,2,\infty\}$.
}

\edef\racor{Corollary \tag}

Let us describe the basic example which will serve as a model for 
the general case.

\proclaim Definition. The Coxeter system $(W,S)$ is of type $B_k$ if its
Dynkin diagram is the following:

$$\pspicture(1.5,0)(9.5,0)
\qline(2,0)(4,0)
\qline(6,0)(8,0)
\qline(8,.05)(9,.05)
\qline(8,-.05)(9,-.05)
\rput(5,0){$\cdots$}
\rput(2.3,.3){$s_k$}
\rput(3.3,.3){$s_{k-1}$}
\rput(4.3,.3){$s_{k-2}$}
\rput(6.3,.3){$s_4$}
\rput(7.3,.3){$s_3$}
\rput(8.3,.3){$s_2$}
\rput(9.3,.3){$s_1$}
\psdots*(2,0)(3,0)(4,0)(6,0)(7,0)(8,0)(9,0)
\endpspicture
\leqno(\mktag)$$

Coxeter group of type $B_k$ is the group of symmetries of a regular
$k$-dimensional cube. Assume that the cube has vertices with each coordinate
equal to $1$ or $-1$. In this case the generator $s_1$ corresponds to the
reflection in the hyperplane defined by vanishing the first coordinate
($s_1$ changes the sign of the first coordinate). The group
$\overline{W_{\{s_1\}}}$ is right-angled Coxeter group generated
by reflections in the hyperplanes defined by vanishing of some coordinate.
We will call such a hyperplane {\it a coordinate hyperplane}.

The parabolic subgroup $W_{\{s_2,\ldots,s_n\}}$ is
the symmetric group of $n$ letters. It acts by permuting coordinates.
Precisely, $s_k$ transposes $(k-1)$st and $k$th coordinates.

\proclaim Proposition \mktag. Let $(W,\{s_1,\ldots,s_m\})$ be of type $B_m$.
If\/ $\Theta$ is a family of $n$ commuting conjugates of $s_1$ then $\Theta$
may be conjugated by an element from $W_{\{s_2,\ldots,s_k\}}$ in such a way
that $\Theta$ consists of {\rm all} conjugates of $s_1$ in
$W_{\{s_1,\ldots,s_n\}}$.

\edef\bpro{Proposition \tag}

{\it Proof: } Family of reflections in coordinate hyperplanes
is defined by the intersection of these hyperplanes. Since all such
intersections can be conjugated by permuting the coordinates,
provided they have the same dimension, the claim follows.\qed


Now return to the general question. Let $\Theta\subset S_T$ be a family
of commuting generators of a right-angled group $\overline{W_T}$.
The aim of the rest of the present  Section is to compute
the f-polynomial of $N_{\overline{W_T}}$ in case when $\overline{W_T}$
is right angled, in terms of the Coxeter matrix of $(W,S)$
(in contrast with \specthm, which gives the f-polynomial in terms of the
growth series of $W$).

In order to compute the f-polynomial of its nerve we need to determine
maximal families of commuting conjugates of generators. 

\proclaim Theorem \mktag. Assume that $T$ is as in \racor.
For any commuting family $\Theta$ of conjugates of generators
there is a unique subset $\Sigma\subset S$ of the same cardinality as $\Theta$,
such that
{\parindent=.3in\parskip=0pt
\item{(1)} $W_\Sigma$ is finite, and if for any $t\in T'=T\cap\Sigma$
the set of vertices of the connected component of
the Dynkin diagram of $W_\Sigma$ containing $t$
is denoted by $\Sigma_t$, then:
\itemitem{$\bullet$} $\Sigma=\bigcup_{t\in T'}\Sigma_t$,
\itemitem{$\bullet$} $\Sigma_t$ is of type $B_{k(t)}$,
\itemitem{$\bullet$} $t$ is the unique element of $\Sigma_t\cap T$ and it is the 
distinguished generator $s_1$;
\item{(2)} there exists $w_\Theta\in W_{S-T}$ such that
$$w_\Theta\Theta w_\Theta=\bigcup_{t\in T'}\{wtw^{-1}|w\in W_\Sigma\}
=\bigcup_{t\in T'}\{wtw^{-1}|w\in W_{\Sigma_t}\}.$$
}

\edef\bbthm{Theorem \tag}

\proclaim Lemma \mktag. If $G$ is a finite subgroup of a Coxeter group,
then it is conjugate to a subgroup of some finite standard parabolic.

\edef\fixlemma{Lemma \tag}

{\it Proof :}
Any Coxeter group acts on its Davis complex [D1] which is CAT(0)
(nonpositively curved) by the theorem of Moussong [M].
The only stabilizers of that action are
the conjugates of finite standard parabolics.

The claim follows, since any action of a finite group on a
CAT(0) complex has a fixed point [BH]. \qed

Obviously, since the set of parabolic subgroups is closed under
the intersection, there exist unique smallest parabolic containing $G$. 

{\it Proof of \bbthm:}
Since a family $\Theta$ of commuting involutions generates
a finite subgroup, by \fixlemma\ 
we may assume that $\Theta$ lies in a finite parabolics $W_{S'}$.
Let $T'$ denote a set of those elements of $T$
whose conjugates appear in $\Theta$.

We claim that $T'\subset S'\cap T$ and moreover that each element
of $\Theta$ is conjugated to some element of $T'$ by an element
of $W_{S'}$. Indeed, by \reflcor\ each element of $\Theta$
(being a reflection) is conjugated inside $W_{S'}$ to some
generator from $S'$. By \racor\ no element of $T$ is a conjugate
of any other generator from $S$, thus the claim.

By \racor\ any generator from $T'$ may be connected in the Dynkin diagram
to any other only by an edge labeled by $4$. Thus by
the classification of finite Coxeter groups ([B]), the connected
component $S'_t$ of $t\in T'$ (within the Dynkin diagram of $W_{S'}$)
is a graph of type $B_k$ for some $k$.

Let $S'_0=S'-\bigcup_{t\in T'}S'_t$. From the definition it follows
that any Coxeter group is the product of its finite parabolics corresponding
to the connected components of its Dynkin diagram.
Thus $W_{S'}=W_{S'_0}\times\bigtimes_{t\in T'}W_{S'_t}$.
Since each element of $\Theta$ is a conjugate of some element
of $T'$, we have $\Theta\subset\bigtimes_{t\in T'}W_{S'_t}$.
Moreover, if $\Theta_t$ denotes those elements from $\Theta$
which are conjugated to $t\in T'$ then $\Theta_t\subset W_{S'_t}$.

By \bpro\ we may find elements $w_t\in W_{S'_t}$
such that each $w_t\Theta_tw_t^{-1}$ consists of all conjugates of $t$
in $W_{\Sigma_t}$ for a suitable subset $\Sigma_t\subset S'_t$.
Since different $W_{S'_t}$ commute, setting $w_\Theta=\prod_{t\in T'}w_t$
completes the proof.
\qed

Note that $\Theta$ determines $\Sigma$ but $\Sigma$ determines $\Theta$
only up to conjugacy. In order to count all commuting families
conjugate to $\Theta$ we need to understand the action of $W_{S-T}$
on the set of families as in \bbthm\ (1).

Let $\Sigma$ be such a family.
The stabilizer of $\Sigma$ is the centralizer of $\Sigma$
in $W_{S-T}$, i.e. $W_{\Sigma_0}$,
where $\Sigma_0=\left(\bigcup_{t\in T}\Sigma_t\right)^\perp-T$.

Summarizing the above we have proved
\proclaim Theorem \mktag. The f-polynomial of\/ $\overline{W_T}$ reads
$$f_{N_{\overline{W_T}},\varphi}=
\sum_{\Sigma\subset S}[W_{S-T}\colon W_{\Sigma_0}]\prod_{t\in T}
{\varphi(t)^{k(t)}\over k(t)!}$$
Where the sum runs over all subsets $\Sigma$ satisfying \bbthm\ (1). 

\edef\develformula{Theorem \tag}

This formula allows the following refinement. Assume that
$\sigma\subset\tau\in N_{\overline W}$.
The corresponding families $\Sigma^\sigma$ and $\Sigma^\tau$
satisfy $\Sigma^\sigma_t\subset\Sigma^\tau_t$ for all $t\in T$.
However if we set $\sigma$ and $\Sigma$ such that
$\Sigma^\sigma_t\subset\Sigma_t$ for all $t\in T$,
then there are several $\tau$ wit the property $\tau\supset\sigma$
and $\Sigma^\tau\subset\Sigma$.

The centralizer of $\Sigma^\sigma$ act on the set of such $\tau$.
thus the number of such $\tau$ equals
$$[W_{\Sigma_0}\colon W_{\Sigma^\sigma_0}].$$

\proclaim Definition. The link ${\rm Lk}_\sigma$
of $\sigma$ a simplex in a simplicial complex $X$
is a subcomplex consisting of all $\tau\in X$ such that
$\tau^*=\sigma\cup\tau\in X$ and $\sigma\cap\tau=\emptyset$.

\proclaim Corollary \mktag. Fix a simplex $\sigma\in N_{\overline{W_T}}$.
It is the commuting family corresponding to $\Sigma^\sigma\subset S$.
The f-polynomial of the link of $\sigma$ reads
$$f_{Lk_\sigma,\varphi}(t)=
\sum_{\Sigma\supset\Sigma^\sigma} [W_{\Sigma^\sigma}\colon W_\Sigma]
\prod_{t\in T}{\varphi(t)^{k(t)-k^\sigma(t)}\over(k(t)-k^\sigma(t))!}.$$

\edef\linkdevelformula{Corollary \tag}

\section An Example.

In this section we would like to present a specific triangulation of
a sphere (which is a boundary complex of a convex polytope),
such that its f-polynomial has non-real roots.
Other such examples are constructed by a different method in [G],
where the reader may also find a connection of the
(former) conjecture on real-rootedness of f-polynomials
of such triangulations to geometry
and combinatorics with further references.

Let $W$ be defined by the following Dynkin diagram:
\def\eightdiag{
\psline[linestyle=dashed](1,2)(2,2)
\qline(2,2)(8,2)
\qline(8,2.05)(9,2.05)
\qline(8,1.95)(9,1.95)
\qline(4,2)(4,1)}
$$\pspicture(1.5,1.75)(9.5,1.75)
\eightdiag
\psdots*(2,2)(3,2)(4,2)(5,2)(6,2)(7,2)(8,2)(4,1)
\psset{dotstyle=o}
\psdots*(1,2)(9,2)
\rput(1.25,2.25){$t_1$}
\rput(9.25,2.25){$t_2$}
\endpspicture
\leqno(\mktag)$$

\edef\eightpic{Diagram \tag}

Let $T=\{t_1,t_2\}$ consist of white dots in the \eightpic.
Recall that $\overline{W_T}=\mathop{\rm ker}(W_S\to W_{S-T}).$

\proclaim Remark \mktag. 
$N_{\overline{W_T}}$ may be realized as a convex triangulation of a sphere,
i.e. a boundary complex of a convex polytope.

{\it Proof :} Consider the group defined by the following Dynkin diagram:
$$
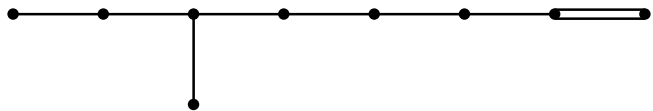
\pspicture(1.5,1.5)(9.5,1.5)
\qline(2,2)(8,2)
\qline(8,2.05)(9,2.05)
\qline(8,1.95)(9,1.95)
\qline(4,2)(4,1)
\psdots*(2,2)(3,2)(4,2)(5,2)(6,2)(7,2)(8,2)(4,1)(9,2)
\endpspicture$$
It can be realized as a finite covolume reflection group in
the hyperbolic space with the fundamental domain $D$ a simplex with
unique ideal vertex [B, Ch. V, \S4, Ex. 17]. Let $\overline D=W_{S_0}\cdot D$.
It is a convex hyperbolic polytope with some ideal vertices.
Finally let $D^*$ be $\overline D$ with ideal vertices truncated.
Is is a polytope dual to $N_{\overline W_T}$.

Since $D^*$ may be realized as a convex polytope in the Klein model
of the hyperbolic space, the polar dual of $D^*$ realizes $N_{\overline W_T}$
as a convex triangulation of a sphere [Z, Cor.~2.14].\qed

According to \develformula\ we have to find all suitable triples
$\{\Sigma_{t_1},\Sigma_{t_2},\Sigma_0\}$.
They are listed in the following picture. Black dots in the connected
component of $t_i$ are elements of $\Sigma_{t_i}$. The remaining black
dots (in the middle) are the element of $\Sigma_0$.

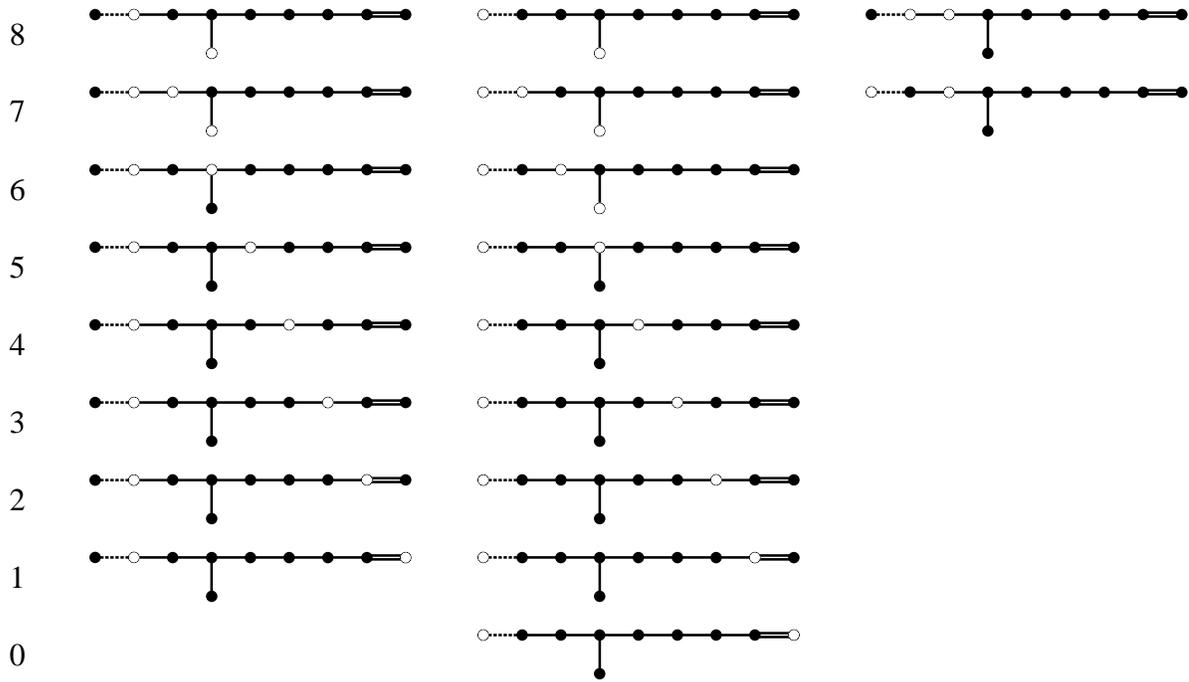
\pspicture(-1.2,0)(8.5,9.3)
\psset{xunit=.43cm,yunit=.43cm,dash=1pt .5pt}
\rput(-1,20){$\sum k(t)$}
\rput(-1,17.5)8
\rput(-1,15.5)7
\rput(-1,13.5)6
\rput(-1,11.5)5
\rput(-1,9.5)4
\rput(-1,7.5)3
\rput(-1,5.5)2
\rput(-1,3.5)1
\rput(-1,1.5)0
\multirput(0,2)(0,2){8}{\eightdiag \psdot(1,2)\psdot[dotstyle=o](2,2)}
\multirput(10,0)(0,2){9}{\eightdiag \psdot[dotstyle=o](1,2)}
\multirput(20,14)(0,2){2}{\eightdiag
\psdots*(4,2)(5,2)(6,2)(7,2)(8,2)(4,1)(9,2)
\psdot[dotstyle=o](3,2)}
\psdots*(12,2)(13,2)(14,2)(15,2)(16,2)(17,2)(18,2)(14,1)
\psdots*(12,4)(13,4)(14,4)(15,4)(16,4)(17,4)(14,3)(19,4)
\psdots*(12,6)(13,6)(14,6)(15,6)(16,6)(18,6)(14,5)(19,6)
\psdots*(12,8)(13,8)(14,8)(15,8)(17,8)(18,8)(14,7)(19,8)
\psdots*(12,10)(13,10)(14,10)(16,10)(17,10)(18,10)(14,9)(19,10)
\psdots*(12,12)(13,12)(15,12)(16,12)(17,12)(18,12)(14,11)(19,12)
\psdots*(12,14)(14,14)(15,14)(16,14)(17,14)(18,14)(19,14)
\psdots*(13,16)(14,16)(15,16)(16,16)(17,16)(18,16)(19,16)
\psdots*(12,18)(13,18)(14,18)(15,18)(16,18)(17,18)(18,18)(19,18)
\psdots*(3,4)(4,4)(5,4)(6,4)(7,4)(8,4)(4,3)
\psdots*(3,6)(4,6)(5,6)(6,6)(7,6)(4,5)(9,6)
\psdots*(3,8)(4,8)(5,8)(6,8)(8,8)(4,7)(9,8)
\psdots*(3,10)(4,10)(5,10)(7,10)(8,10)(4,9)(9,10)
\psdots*(3,12)(4,12)(6,12)(7,12)(8,12)(4,11)(9,12)
\psdots*(3,14)(5,14)(6,14)(7,14)(8,14)(4,13)(9,14)
\psdots*(4,16)(5,16)(6,16)(7,16)(8,16)(9,16)
\psdots*(3,18)(4,18)(5,18)(6,18)(7,18)(8,18)(9,18)
\psdots*(21,18)(22,16)
\psset{dotstyle=o}
\psdots*(19,2)(18,4)(17,6)(16,8)(15,10)(14,12)(13,14)(14,13)
(12,16)(14,15)(14,17)
\psdots*(9,4)(8,6)(7,8)(6,10)(5,12)(4,14)(3,16)(4,15)(4,17)
\psdots*(21,16)(22,18)
\endpspicture

The orders of finite Coxeter groups are well known ([B])
Substituting the above to \develformula\ one obtains
$$\eqalign{f_{N_{\overline{W_T}}}(t_1,t_2)
&=696729600\big(({2t_2^7t_1\over5040}+{t_2^8\over40320})+
({t_2^6t_1\over720}+{t_2^7\over5040}+{t_2^7\over2\cdot5040})\cr
&\hphantom{=696729600\big(}+({t_2^5t_1\over2\cdot2\cdot120}
+{t_2^6\over2\cdot720})+
({t_2^4t_1\over24\cdot24}+{t_2^5\over2\cdot6\cdot120})\cr
&\hphantom{=696729600\big(}+({t_2^3t_1\over6\cdot192}+{t_2^4\over24\cdot120})+
({t_2^2t_1\over2\cdot1920}+{t_2^3\over6\cdot1920})\cr
&\hphantom{=696729600\big(}+({t_2t_1\over23040}+{t_2^2\over2\cdot51840})
+({t_1\over322560}+{t_2\over2903040})
+{1\over696729600}\big)\cr
&=t_1(276480\,t_2^7+967680\,t_2^6+1451520\,t_2^5+1209600\,t_2^4+604800\,t_2^3\cr
&\quad+181440\,t_2^2+30240\,t_2+2160)\quad+(17280\,t_2^8+207360\,t_2^7\cr
&\quad+483840\,t_2^6+483840\,t_2^5+241920\,t_2^4+60480\,t_2^3+6720\,t_2^2
+240\,t_2+1).\cr}$$
Therefore, in particular,
$$\eqalign{f_{N_{\overline{W_T}}}(t,t)
&=293760\,t^8+1175040\,t^7+1935360\,t^6+1693440\,t^5+846720\,t^4\cr
&\quad+241920\,t^3+36960\,t^2+2400\,t+1.\cr}$$
and
$${1\over\overline{W_T}(t)}={
\displaystyle{t^8-2392*t^7+20188*t^6-70504*t^5+107590*t^4\qquad\atop
\hfill-70504*t^3+20188*t^2-2392*t+1
}
\over (1+t)^8}
$$
Thus the poles of $\overline{W_T}(\cdot)$
(with 2 digit precision) are: $0.41\cdot10^{-3}$,
$0.24\pm0.16i$, $0.63$, $1.6$, $2.9\pm1.9i$, $2.4\cdot10^3$.
This provides a first known counterexample to the Real
Roots Conjecture [G].

The link $L$ (of codimension 2) of a face
$$
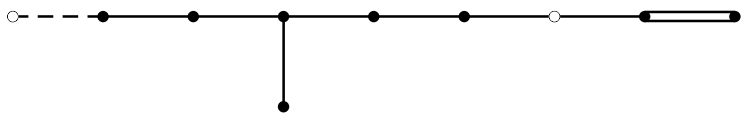
\pspicture(1.5,1)(9.5,2)
\eightdiag
\psdots*(2,2)(3,2)(4,2)(5,2)(6,2)(8,2)(9,2)(4,1)
\psset{dotstyle=o}
\psdots*(1,2)(7,2)
\endpspicture$$
is a triangulation of a five dimensional flag sphere.
Using \linkdevelformula\ we check that $f_L(-1/2)\neq 0$.
A link $K$ (in $L$) of an edge defined by
$$
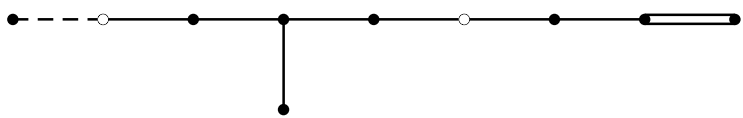
\pspicture(1.5,1)(9.5,2)
\eightdiag
\psdots*(1,2)(3,2)(4,2)(5,2)(7,2)(8,2)(9,2)(4,1)
\psset{dotstyle=o}
\psdots*(2,2)(6,2)
\endpspicture$$
satisfies $f_K(t)=(1+2t)^4$, which is an f-polynomial
of a cross polytope (and, as follows, $K$ is a cross polytope).

Such a pair $(K,L)$ allows us to construct a triangulation of a
five dimensional sphere which is a counterexample to the
Real Roots Conjecture [G, Th.~5.2].
Five is also the smallest dimension when this is possible [G, Th.~3.2].

\section References

\parindent=.4in
\item{[B]} N.~Bourbaki, {\it Groupes et alg\`ebres de Lie, chapitres IV-VI},
Hermann 1968,
\item{[BH]} M.~R.~Bridson and A.~Haefliger, {\it Metric spaces of
non-positive curvature}, Springer, Berlin, 1999.
\item{[D1]} M.~W.~Davis, {\it Nonpositive curvature and reflection groups},
Handbook of Geometric Topology, R.~Daverman and R.~Sher ed., Elsevier 2001.
\item{[D2]} M.~W.~Davis, book in preparation
\item{[G]} S.~R.~Gal, {\it Real Root Conjecture fails for five and higher dimensional spheres},
to appear in Discrete \& Computational Geometry
\item{[M]} G.~Moussong, {\it Hyperbolic Coxeter groups}, PhD Dissertation,
Ohio State University, 1987.
\item{[S]} J.~P.~Serre, {\it Cohomologie des groupes discrets}, in
{\it Prospects in Mathematics}, pp.~77-169,
Annal of Math.~Studies No. {\bf 70}, Princeton 1971
\item{[T]} J.~Tits, {\it Le probl\`eme des mots dans les
groups de Coxeter} in Symposia Mathematica, vol.~1, London 1969
\item{[Z]} G.~M.~Ziegler, {\it Lectures on polytopes},
Graduate Texts in Mathematics, Vol.~152, Springer-Verlag, New York 1995

\bigskip

\rightline{\it Wroc\l aw, April 2002}
\bye